\documentclass[conference]{IEEEtran}
\newcounter{MYtempeqncnt}
\usepackage{slashbox}
\usepackage{cite}
\usepackage{amssymb}
\ifCLASSINFOpdf
   \usepackage[pdftex]{graphicx}

\usepackage{amsmath,bm}

\usepackage{array}

\begin{document}
\title{An LMI Based Stability Margin Analysis for Active PV Power Control of Distribution Networks with Time-Invariant Delays}
\author{\IEEEauthorblockN{Chang Fu, Caisheng Wang, Saeed Alyami}
\IEEEauthorblockA{Department of Electrical and Computer Engineering,
Wayne State University, Detroit, Michigan USA\\
 chang.fu@wayne.edu; cwang@wayne.edu; saeed.alyami2@wayne.edu}}

\maketitle

\begin{abstract}
 High penetration of photovoltaic (PV) generators can lead to voltage issues in distribution networks. Various approaches including the real power control through PV inverters have been proposed to address voltage issues. However, among different control strategies, communication delays are inevitably involved and they need to be carefully considered in the control loop. Those delays can significantly deteriorate the system performance with undesired voltage quality, and may also cause system instability. In this paper, according to the inverter based active power control strategy, a linearized state space model with communication delay is presented. A delay dependent stability criterion using linear matrix inequality (LMI) approach is used to rigorously obtain the delay margins based on different system parameters. The method can handle multiple PVs in the distribution network as well.
\end{abstract}


\IEEEpeerreviewmaketitle

\section{Introduction}
\IEEEPARstart{A}{s} one of the most important clean renewable sources for sustainable energy development, PV generation has been rapidly increased for more than two decades worldwide \cite{16}. 

 The majority of the PV systems has been and will be installed in distribution networks. As a result, the PV penetration level will become unprecedentedly high (e.g. well over 50\%) and continue to grow around the world \cite{41}. The high penetration of PV systems has led to great technical challenges, including voltage problems, harmonics, grid protection, etc., in the operation and development of modern distribution networks \cite{20,21}. \par

A poor voltage profile in a distribution network may lead to issues of power quality, equipment safety, system reliability and stability, and thus can raise system losses and cause equipment damages. Overvoltage is one of the most significant concerns among the above mentioned challenges, and limits the capacity of PV accommodation. Overvoltage issues may happen when the solar irradiance is high while the load demand is low so that the voltages at certain nodes may exceed the upper acceptable limit due to reverse power flows. 
 Active power curtailment (APC) methods have  been widely studied to address overvoltage issues by exploring the real power control capability of PV inverters, such as droop based active power curtailment, global voltage sensitivity matrix method, adaptive real power capping method, as well as consensus based method \cite{01,13,27,28,29,30}. 
However, most of those approaches are based on the assumption that the control signals and all the measurements are obtained, processed and delivered in an ideal communication environment with no communication delays.  Due to large numbers of components in distribution networks, it is not economically feasible to have dedicated channels for communications among local control devices and between local control devices and the central controller. For communications over a shared channel, delays and losses are inevitable, which introduce a great challenge for control of distribution networks with high penetration of PVs.\par

Linear matrix inequality (LMI) based stability methods have been studied extensively in recent decades. According to \cite{36,37,38}, different stability criteria with respect to different types of communication delays, such as time-invariant delay and time-variant delay, were investigated.  In \cite{39,40}, LMI based state feedback controller designs were provided according to the state space mode,l and a controller $K$ was designed to stabilize the system with communication delays. In power systems, LMI approach has been used in load frequency control to mitigate the impact of communication delays existing at the area control error (ACE) signals. A frequency regulation controller was designed in \cite{4} based on an asymptotically stable LMI constraint illustrated in \cite{5}. An LMI based design was also used in \cite{6} to derive different controller parameters for frequency regulation of large power systems under different communication delays. \par
For distribution networks with high penetration of PVs, it is necessary to have a controller that can guarantee the systems' stability and performance while the systems are subject to various disturbances (load and solar irradiance variations) and communication delays. Therefore, in this paper, in addition to control the voltage in the distributed network via active power, an LMI based stability criterion will also be studied according to the state space model considering communication delays. \par

The rest of the paper is organized as follows: The system under study is discussed and the model of multiple PV connected distribution network with  communication delay is given in detail in Section \ref{modeldetail}. The LMI based delay-dependent time-invariant stability criterion is given in Section \ref{sec3}. Simulation studies are carried out based on the active power control method in a small distribution system in Section \ref{casestudy}, followed by the calculations of delay margins with respect to different system parameters. The conclusion is drawn in Section \ref{conclusion}.

\section{ Model of a grid-tied PV system with power control capability}
\label{modeldetail}
The newly developed smart PV inverters are equipped with the capability of regulating real power between zero and the maximum power point (MPP) and providing reactive power compensation \cite{35}. The PV systems can then operate as controllable sources. Nevertheless, most of the distributed PV systems have limited reactive power regulation capabilities such as adjusting power factor between 0.95-lead and 0.95-lag due to economic considerations \cite{35}. Moreover, the focus of this paper is on the impact analysis and mitigation of communication delays for voltage regulation of distributed PVs. Therefore, only real power control is considered for PV systems in this paper. Reactive power regulation with communication delays can be studied in a similar way.\par
\begin{figure}[!t]
\centering
\begin{minipage}[t]{1\linewidth}
    \includegraphics[width=\linewidth,height=1.9in]{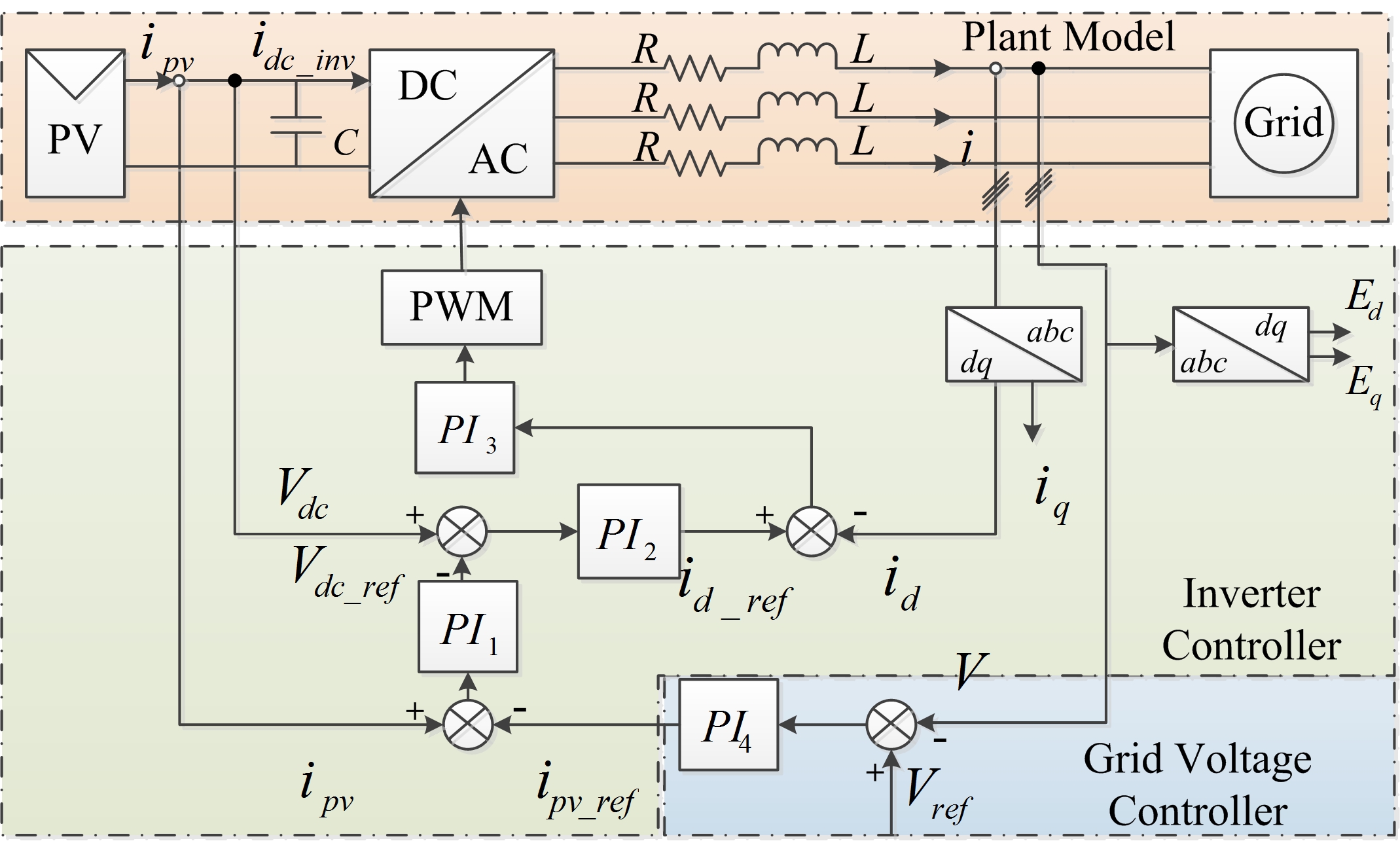}
\caption{Schematic diagram of a grid-connected PV system. }
\label{detailschematic}
\end{minipage}
\end{figure}
The schematic diagram of a grid-tied PV system with real power control capability for voltage regulation is shown in Fig. 1. The system has three major parts: the plant model, the inverter controller and the grid voltage controller. The plant model is the main power circuit of the system. There is no DC/DC converter in the circuit and the DC/AC inverter is used to not only convert DC power into AC, but also achieve real power control via the inverter controller. The inverter controller regulates the PV output current ($i_{pv}$) to follow the PV current reference ($i_{pv\_ref}$) by generating a dc bus voltage reference ($V_{dc\_ref}$). Since there is no DC/DC converter, the dc bus voltage $V_{dc}$ is also the PV output voltage. The dc bus voltage error signal ($V_{dc}-V_{dc\_ref}$) is then used to regulate the inverter output real power represented by the d-axis current $i_d$.  The grid voltage controller takes the grid voltage control error ($V_{ref}-V$) as input and generates the reference current ($i_{pv\_ref}$) for controlling the inverter. For the purpose of analysis, four PI controllers (represented by $PI_i$ in Fig. \ref{detailschematic}) are used in this model though some other types of controller can be used as well. $ PI_1$ is used to produce the reference dc voltage of the PV panel $V_{dc\_ref}$ from ($i_{pv} -i_{pv\_ref}$ ). $PI_2$ and $PI_3$ are traditional controllers usually used in the DC/AC inverter control. $PI_4$ generates the reference current $i_{pv\_ref}$ of the PV panel from ($V_{ref}- V$). The state variables of the system are chosen as follows:\par
$x_1$: the integrator output of $PI_1$\par 
$x_2$: the integrator output of $PI_2$\par
$x_3$: the first order response of $i_{d\_ref}$\par
$x_4$: the integrator output of $PI_3$\par
$x_5$: the d-axis component of the output current $i_d$\par
$x_6$: the input current of the inverter $i_{dc\_inv}$\par
$x_7$: the voltage of the dc side of the inverter $V_{dc}$, $V_{dc}=V_{pv}$\par
$x_8$: the integrator output of $PI_4$\par
\begin{figure*}[!t]
\begin{minipage}{0.68\textwidth}
\includegraphics[width=0.95\linewidth,,height=2.0in]{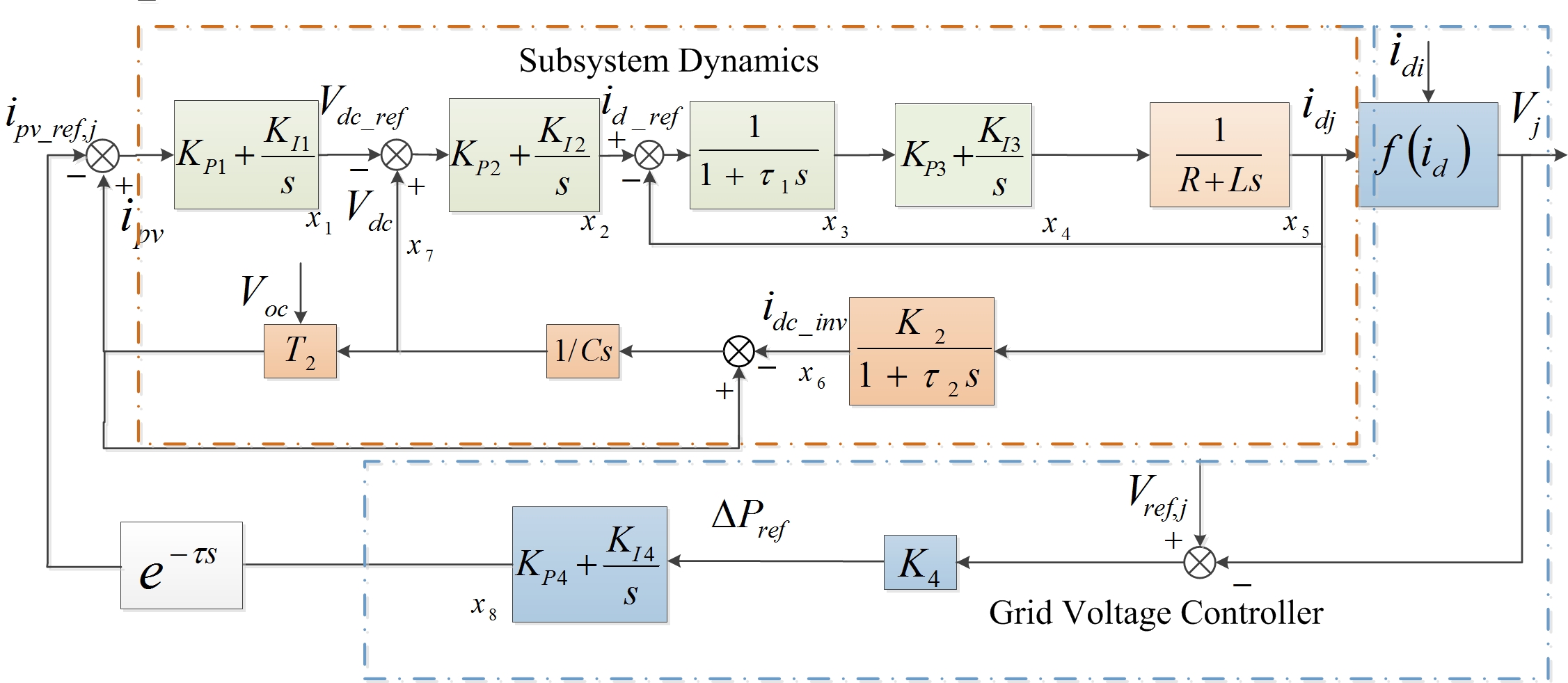}
\caption{Detailed model of the $j$th grid-connected PV system.}
\label{systemmodel3}
\end{minipage}\hfill
\begin{minipage}{0.32\textwidth}
\includegraphics[width=1\linewidth,,height=2.0in]{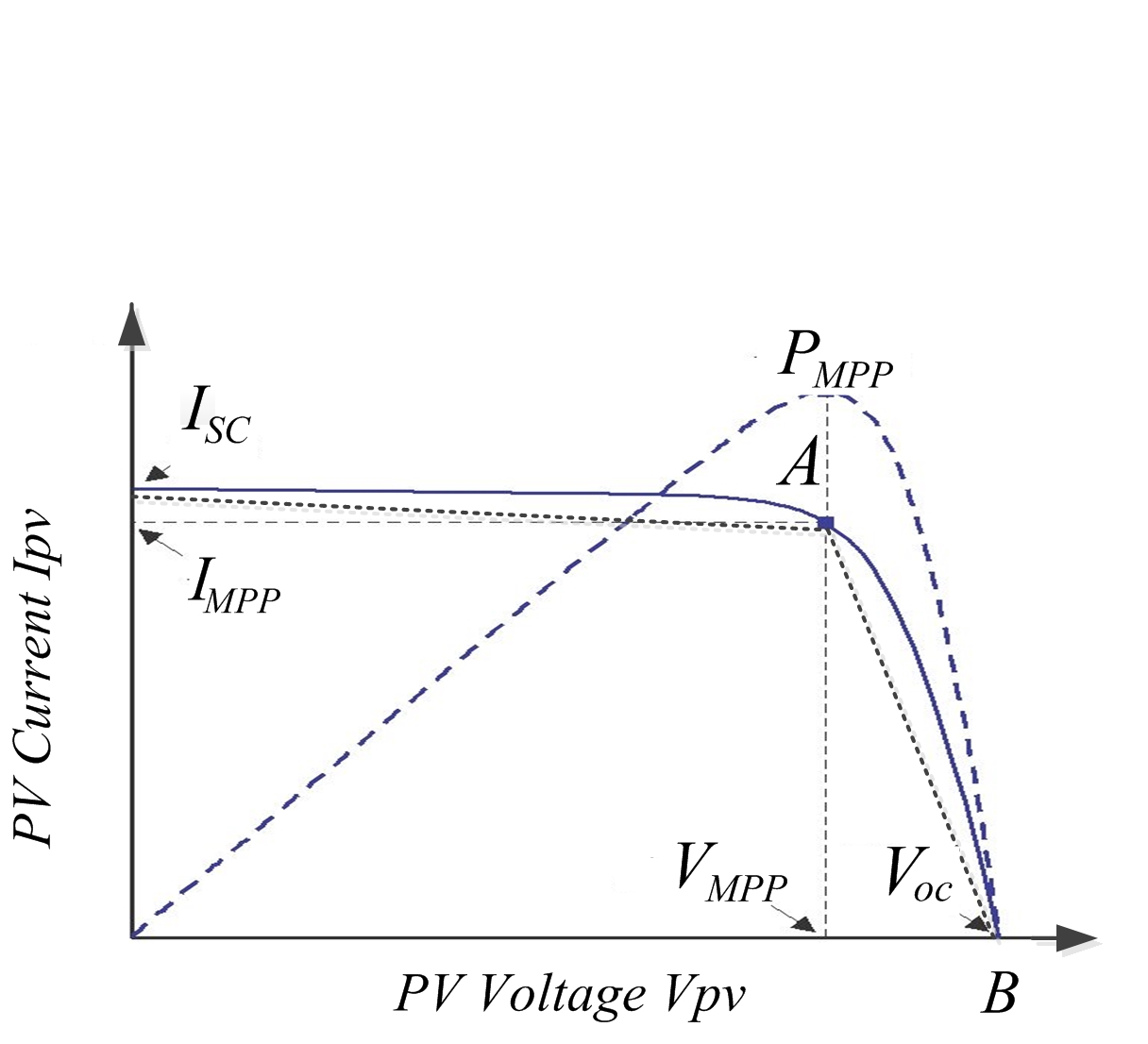}
\caption{PV I-V curve and the linearization of $R_{pv}$.}
\label{pvcurve}
\end{minipage}
\end{figure*}
The detailed dynamic model of the PV system is shown in Fig. \ref{systemmodel3}. For analysis, the system is a linearized model. It should be noted that the model is for the $j$th PV system. The subscript $j$ is omitted for the cases where the omission will not cause confusion. In the figure, the four PI controllers are represented by ($K_{Pi} +K_{Ii}/s$), $i=1,...,4$. $T_2$ is a linearized model for describing the $V-I$ characteristic of the PV panel. The PV system normally works between $V_{MPP}$ and $V_{oc}$, as shown in Fig. \ref{pvcurve}. $V_{oc}$ is the open circuit voltage of the $PV$ panel and $V_{MPP}$ is the voltage at the maximum power point. $V_{oc}$ is a function of solar irradiance, temperature, the material of the panel, the number of PVs connected in series, etc. For a fixed PV panel,  if the solar irradiance and the temperature are constant, then $V_{oc}$ is also a constant. The $V-I$ relationship of the PV panel can be represented as:
\begin{equation}
i_{pv}=(V_{oc}-V_{pv})/R_{pv}
\end{equation}
where $V_{pv} =V_{dc}$. $R_{pv}$ is the equivalent resistance of the panel and changes with PV voltages. For the purpose of simplification without introducing a large error, a linearized equivalent resistance $R_{pv}$ is used in this model. In other words, a straight line $AB$ is considered in calculating $R_{pv}$, shown in Figs. \ref{systemmodel3} and \ref{pvcurve}.\par
\begin{equation}
R_{pv}=\frac{V_{oc}-V_{MPP}}{I_{MPP}}
\label{rpvformula}
\end{equation}
where $I_{MPP}$ is the corresponding PV current when the PV works at its maximum power point.\par 
The dc input current to the inverter is determined by the inverter output current, and ultimately by the inverter output power. The relationship between the input and output currents of the inverter is modeled by a first order system in this study.\par
\begin{equation}
i_{dc\_inv}=i_d\frac{K_2}{1+\tau_2s}
\label{currentdynamic}
\end{equation}
 where $i_d$ is the d-axis component of the inverter output current, $i_{dc\_inv}$ is the dc input current of the inverter, and $\tau_2$ is the time constant of the first order system. $K_2$ is a constant, which can be calculated based on power balancing between the input and output powers of the inverter. Assuming there is no power loss in the inverter and the d-axis is aligned with phase A, the following power balance equation can be established when there is no reactive power (i.e., $ iq = 0$, $E_q=0$):
\begin{equation}
\nonumber
V=\sqrt{E^2_d+E^2_q}=E_d
\end{equation}
\begin{equation}
V_{dc}\cdot i_{dc\_inv}=\frac{3}{2}E_di_d
\Longrightarrow V_{dc\_ref}\cdot i_{dc\_inv}=\frac{3}{2}E_di_d
\end{equation}
where $E_d$ is the d-axis component of the grid voltage, and $E_q$ is the q-axis component of the grid voltage.\par
In this case, $K_2$ can be obtained as:
\begin{equation}
K_2=\frac{3E_d}{2V_{dc\_ref}}
\end{equation}\par
The block of $f(\bm{i_d})$, $\bm{i_d}=[i_{d1},...,i_{dn}]^T$,  in Fig. \ref{systemmodel3} is used to represent the power network algebraic equations that link the injected currents and the system bus voltages. The system bus voltages ($\boldsymbol{V_{bus}$}) can be obtained based on the system impedance matrix  ($\boldsymbol{Z_{bus}$})  and the injected current sources  ($\boldsymbol{I_{bus}$}) as \cite{book1}:
\begin{equation}
\bm{V_{bus}=Z_{bus}I_{bus}}
\label{nodevoltage}
\end{equation} \par
For instance, the voltage of bus $j$  in the power network can be obtained as:
\begin{equation}
\label{volt}
V_j=\sum_{i=1}^nZ_{ij}i_{di}+\sum_{i=n+1}^NZ_{ij}i_{k}=\sum_{i=1}^nZ_{ij}i_{di}+V_j^0
\end{equation}
where $i_{di}$ ($i=1,...,n$) are the injected currents from the $PV$ sources and $n$ is the total number of PVs in the distribution grid. $V_j^0$ is the voltage contribution to bus $j$ from other generation sources $i_k$, which is a constant in this study.\par
$K_4$ is a constant coefficient that converts the voltage error ($V_{ref,j} - V_j$) into a power change signal to control the $PV$ system. From the control viewpoint. In practice, $K_4$ can be obtained from the Jacobian matrix of the network.
Therefore, based on the above assumptions, linearizations and simplifications, the state space model of the proposed system can be developed, which is composed by the PV subsystem dynamics and the grid voltage controller:
\subsubsection{Subsystem Dynamics}
The subsystem dynamics describe the transient characteristics of the PV systems and the inverter controller. As shown in Fig. \ref{systemmodel3}, the state space representation of the subsystem dynamics of PV$_j$ can be written as:\par
\begin{equation}
\label{subsysdyn}
\begin{cases}
\dot{x_j}(t)=A_jx_j(t)+B_ju_j(t)+H_j\sigma_j(t)\\
y_j(t)=C_jx_j(t)+D_ju_j(t)
\end{cases}
\end{equation}
The $u_j$ is the input to the subsystem, and $i_{pv\_ref,j}$ is the output of the grid voltage controller. The output $y_j$ is the d-axis component of the PV output current $i_{dj}$, and $\sigma_j$ is the disturbance. The disturbance is represented as the open circuit voltage ($V_{oc}$ in Fig. \ref{systemmodel3}) of the PV panel in this study. $V_{oc}$ will change as the solar irradiance varies.

\subsubsection{Grid Voltage Controller}
 The state space  representation of the grid voltage controller is shown in (\ref{localcontroller}):
\begin{equation}
\label{localcontroller}
\begin{cases}
\dot{z_j}(t)=A{^F_j}z_j(t)+B{^F_j}\omega_j(t)\\
u_j(t)=C{^F_j}z_j(t)+D{^F_j}\omega_j(t)
\end{cases}
\end{equation}
where $\omega_j$ is the input of the local controller and it equals to ($V_{ref,j}-V_j$). $u_j$ is the output of the controller, which is $i_{pv\_ref,j}$, and it is used as the input in (\ref{subsysdyn}). $A{^F_j}$, $B{^F_j}$, $C{^F_j}$ and $D{^F_j}$ can be found in (\ref{af}):
\begin{eqnarray}
A{^F_j}=\begin{bmatrix}
0\\
\end{bmatrix}
B{^F_j}=\begin{bmatrix}
\label{af}
K_4K_{I4}\\
\end{bmatrix}
C{^F_j}=\begin{bmatrix}
1
\end{bmatrix}
D{^F_j}=\begin{bmatrix}
K_4K_{P4}\\
\end{bmatrix}
\end{eqnarray}\par

In a distribution network with a central controller, the other PVs in the network may be subject to similar communication delays when there is a control signal sent from the same controller. A single delay $e^{-\tau s}$ is used in this study to approximate such scenarios. Therefore, the input to the subsystem shown in Fig. \ref{systemmodel3} can be written as $u_j(t-\tau)$. Take the 2-PV system in Fig. \ref{disnet} as the example. Substitute (\ref{localcontroller}) and (\ref{volt}) into (\ref{subsysdyn}), the state space model of the delayed system can be written as:
\begin{equation}
\label{ssdelay}
\begin{cases}
\dot{x^{\prime}}(t)=Ax^{\prime}(t)+A_dx^{\prime}(t-\tau)+F\sigma^{\prime}(t)\\
y^{\prime}(t)=Cx^{\prime}(t)
\end{cases}
\end{equation}
where 
\begin{eqnarray}
\nonumber
A&=\begin{bmatrix}
A_1&0&0&0\\
0&A_2&0&0\\
-B_1^FZ_{66}C_1&-B_1^FZ_{67}C_2&A_1^F&0\\
-B_2^FZ_{76}C_1&-B_2^FZ_{77}C_2&0&A_2^F
\end{bmatrix}\\
\nonumber
\end{eqnarray}
\begin{eqnarray}
\nonumber
A_d&=&\begin{bmatrix}
-B_1D_1^FC_1Z_{66}&-B_1D_1^FC_2Z_{67}&B_1C_1^F&0\\
-B_2D_2^FC_1Z_{76}&-B_2D_2^FC_2Z_{77}&0&B_2C_2^F\\
0&0&0&0\\
0&0&0&0
\end{bmatrix}\\
\nonumber
F&=&\begin{bmatrix}
H_1&B_1D_1^F&0&0\\
0&B_1^F&0&0\\
0&0&H_2&B_2D_2^F\\
0&0&0&B_2^F
\end{bmatrix}
\nonumber
C=diag[C_1,C_2,0,0]
\end{eqnarray}
$Z_{ij}$ can be found in (\ref{nodevoltage}), new state variable $x^{\prime}(t)$ is a combination of $x_j(t)$ and $z_j(t)$, $\sigma^{\prime}(t)$ is the new system disturbance composed by $\sigma_j(t)$ and $V_{ref,j}(t)$. Other notations are similar as those in (\ref{localcontroller}) and (\ref{subsysdyn}). For system with $n$ PVs installed, the state space model can be extended by the same approach.
\section{Delay-dependent Time-invariant Stability Criterion}
\label{sec3}
 
Consider a system with time delay
\begin{equation}
\setcounter{MYtempeqncnt}{\value{equation}}
\setcounter{equation}{12}
\dot{x}(t)=Ax(t)+A_dx(t-\tau), \tau\geq 0
\end{equation}
where $\tau$ is the time delay. The system stability holds for $\tau<\tau_d$, where $\tau_d$ is the stability margin, and for $\tau>\tau_d$, the system is unstable. Many methods can be used to calculate $\tau_d$. A delay-dependent time-invariant stability criterion proposed in \cite{36} can be used to determine the delay margin of a distribution network with PVs installed:\par
$\bm{Theorem\ 1}$: Assume that an uncertain time-invariant time delay in [0, $\tau_d$], i.e., $\tau\in[0,\tau_d]$. Then if there exists $P>0$, $Q>0$, $V>0$ and $W$ such that 
\begin{equation}
\label{lmi}
\begin{bmatrix}
(1,1)&-W^TA_d&A^TA_d^TV&(1,4)\\
-A_d^TW&-Q&A_d^TA_d^TV&0\\
VA_dA&VA_dA_d&-V&0\\
(1,4)^T&0&0&-V
\end{bmatrix}<0
\end{equation}
where
\begin{equation}
\nonumber
(1,1)\triangleq(A+A_d)^TP+P(A+A_d)+W^TA_d+A_d^TW+Q
\end{equation}
\begin{equation}
\nonumber
(1,4)\triangleq\tau_d[W^T+P]
\end{equation}
then the system is asymptotically stable. The proof of this theorem can be found in \cite{36}. 
\section{Case Study}
\label{casestudy}
In this section, simulation studies are carried out based on the proposed active power control method. The LMI based stability criterion is also studied to calculate the delay margin of the PV connected distributed system. The effectiveness of the voltage regulation method is verified in a delay free system, and the delay margins are calculated according to (\ref{ssdelay}) and (\ref{lmi}), by using different system parameters.  
\subsection{Active Power Voltage Regulation}
\begin{figure}[!t]
\begin{minipage}[t]{1\linewidth}
    \includegraphics[width=\linewidth]{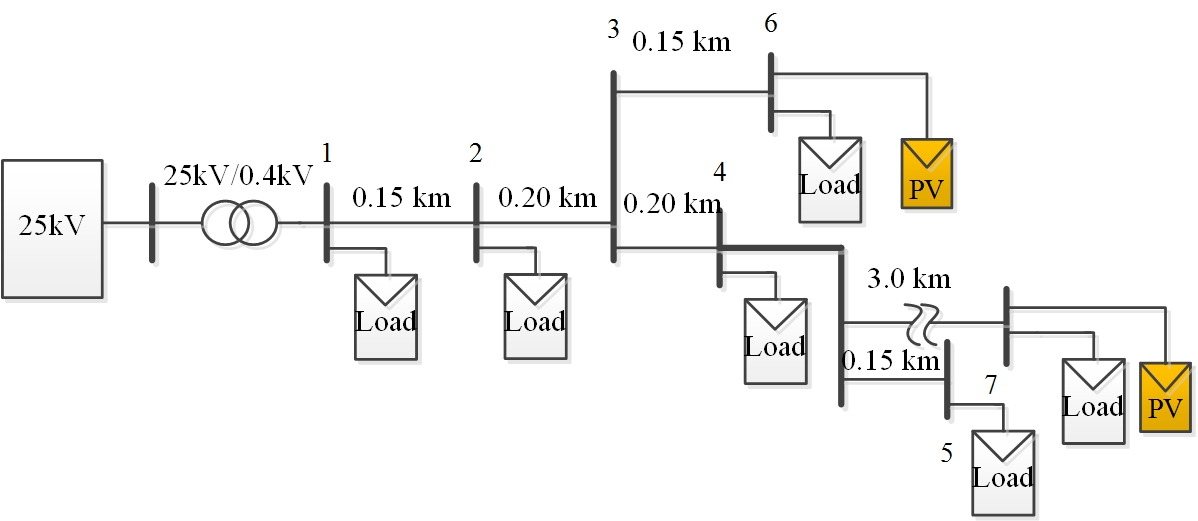}
\caption{Tested distribution network with six sub-communities. }
\label{disnet}
\end{minipage}
\end{figure}

The simulation study is carried out in a small test network in Fig. \ref{disnet}. The test network is a residential suburban feeder with 6 sub-communities. The LV feeder is connected to the grid through a 25kV/0.4kV transformer. Each load in Fig. \ref{disnet} represents an individual sub-community with a different level of power consumption. The PV systems are connected to Nodes 6 and 7 to provide additional power to the grid. Nodes 1 to 6 are relatively close to each other while Node 7 is connected to the network via a 3.0 km power line. The solar irradiance is set to $1000W/m^2$ at 25$^{\circ}$C to simulate a condition that may generate an undesired voltage profile. The sizes of the loads are also shown in Fig. \ref{disnet}.\par

\begin{figure}[!t]
\begin{minipage}[t]{1\linewidth}
    \includegraphics[width=\linewidth]{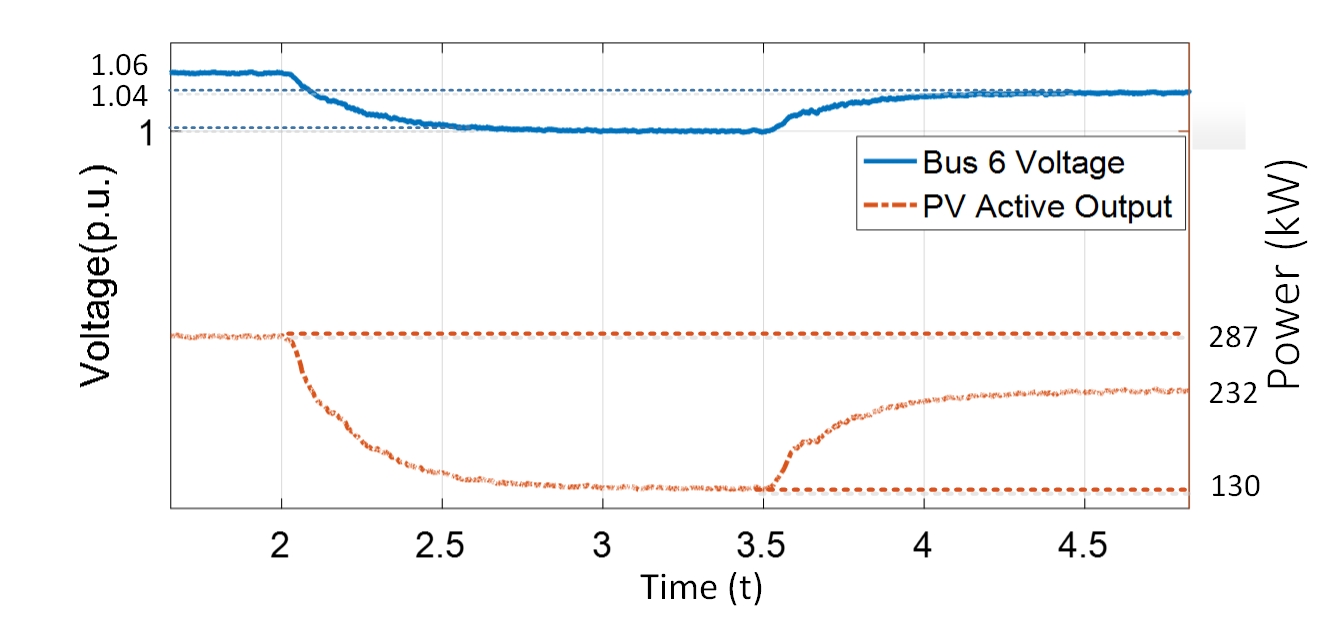}
\caption{Voltage and PV active power output at Node 6. }
\label{voltage}
\end{minipage}
\end{figure}

The performance of the PV active power curtailment method is shown in Fig. \ref{voltage}, where the voltage profile of Node 6 is 1.06 p.u. (solid line), which already exceeds the critical value (1.05 p.u.) set for this study. Overvoltage may cause damage to the electrical components and the proposed active power control method is applied to control the PV inverter to regulate the voltage back under the critical value. At $t$=2s, a reference voltage $V_{ref}$=1.00 p.u. is sent from the center controller, and the PV at Node 6 curtails its output active power from 287 kW to 130 kW, and the voltage at Node 6 successfully reaches the reference value at $t$=2.6s. At $t$=3.5s, a new reference value of $V_{ref}$=1.04 p.u. is set. The PV output then raises to 230 kW, and the voltage of Node 6 shown in Fig. \ref{voltage} is also increased to 1.04 p.u..
\subsection{Delay-dependent Stability Criterion}
\begin{figure}[!t]
\begin{minipage}[t]{1\linewidth}
    \includegraphics[width=\linewidth]{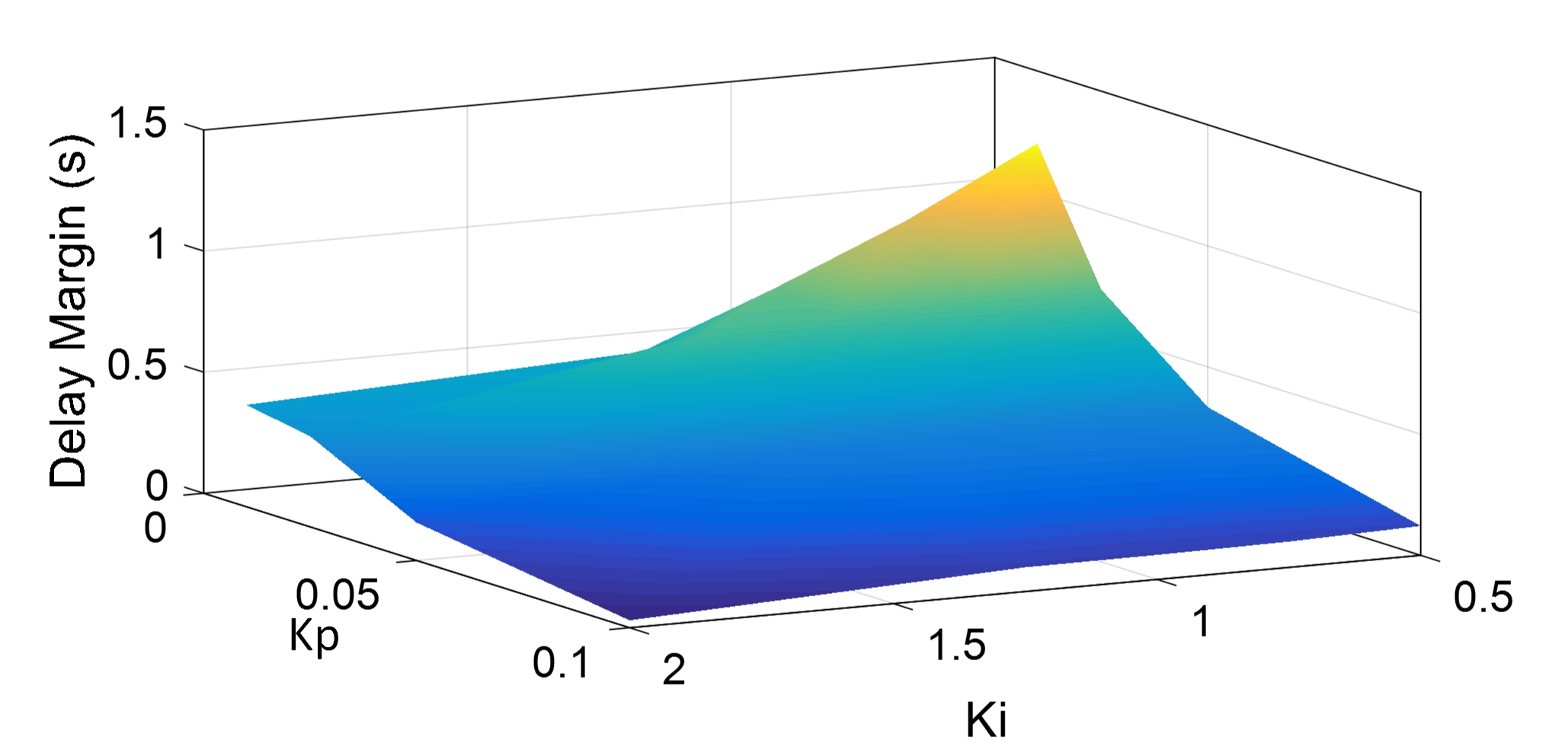}
\caption{Relationships among $\tau_d$, $K_{P4}$ $\&$ $K_{I4}$. }
\label{voltage2}
\end{minipage}
\end{figure}

The delay-dependent time-invariant stability criterion can be obtained by solving the LMI in (\ref{lmi}) with different sets of gains of the PI controllers. Table \ref{delaymargin} shows the different delay margins with  respect to different $K_{P4}$ and $K_{I4}$. The results of $\tau_d$ indicate that for a constant communication delay, the delay margin $\tau_d$ increases with the decreasing of $K_{P4}$ and $K_{I4}$, especially when $K_{P4}$ and $K_{I4}$ are small (e.g., $K_{P4}$=0.01, $K_{I4}=0.5$). As shown in Fig. \ref{voltage2}, a very sharp increase can be found when the corresponding parameters are relatively small.
\begin{table}[!h]
\caption{Different Delay Margins with respect to $K_{P4}$ $\&$ $K_{I4}$ }  
\label{delaymargin}
  \centering
\begin{tabular}{| c|c|c|c|c|c |}
\hline
\backslashbox{$K_{P4}$}{$K_{I4}$}&0.5&0.75&1.00&1.25&2\\
  \hline
\hline
 0.100 &0.121&0.106&0.101&0.099&0.030\\
0.050 & 0.334 &0.315&0.305&0.287&0.156 \\
 0.025& 0.681&0.670&0.652&0.447&0.374 \\
0.010& 1.200&0.931&0.673&0.491&0.423\\ 
\hline 
\end{tabular}
\end{table}\par
\begin{figure}[!t]
\begin{minipage}[t]{1\linewidth}
    \includegraphics[width=\linewidth]{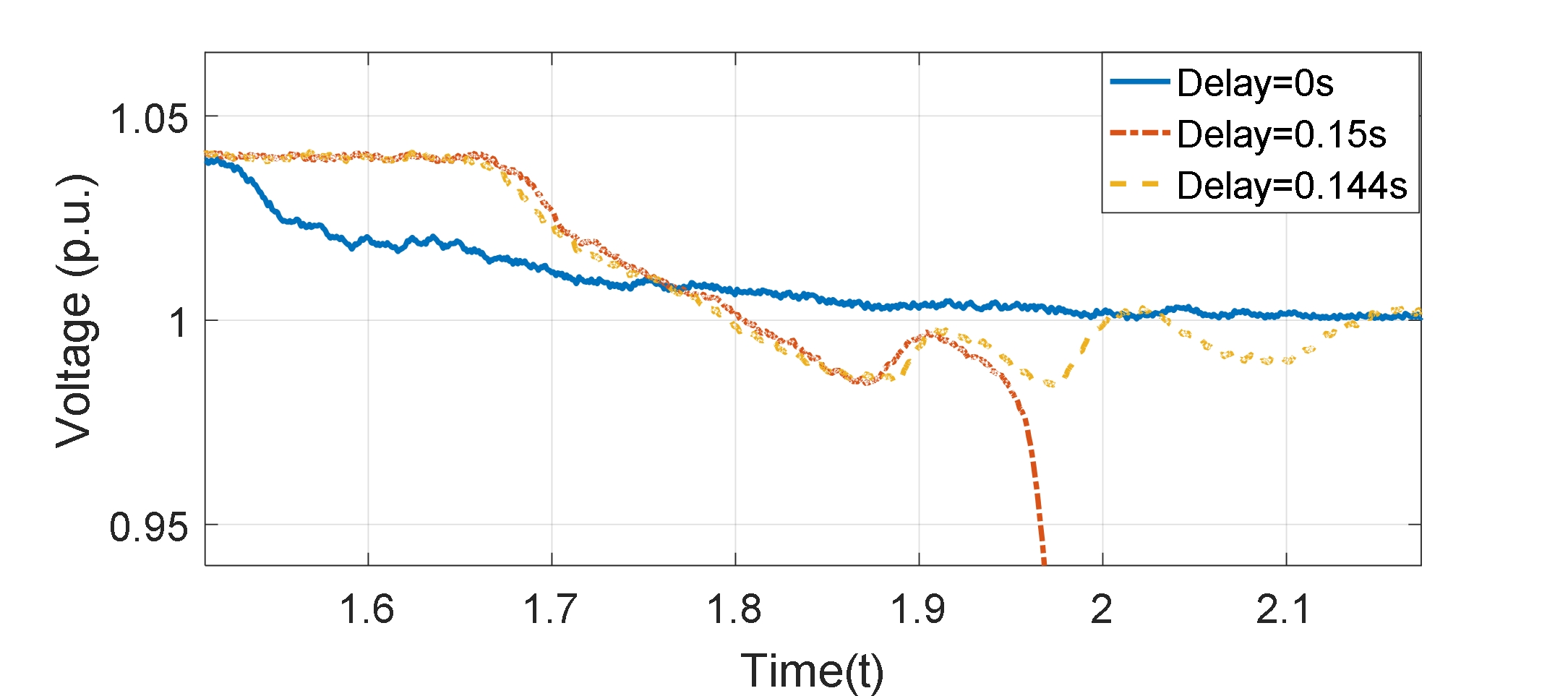}
\caption{Voltage profile at Node 6 with different communication delays.}
\label{compare}
\end{minipage}
\end{figure}

The simulation study has also been carried out to verify the accuracy of the calculated delay margin according to the linearized model. Due to the linearization,  there is a small error between the calculated value $\tau_d$ and the real value obtained from the simulation study. For instance,  the delay margin is calculated as 0.156s when $K_{P4}$=0.05, $K_{I4}$=2 (shown in Table \ref{delaymargin}).  According to the simulation result given in Fig. \ref{compare}, the real delay margin is found to be 0.144s. When the delay exceed the cirtical value (e.g., $\tau$=0.15s), the voltage profile shown in Fig. \ref{compare} (dash-dot) indicates the system becomes unstable. The delay margin obtained can help set the upper bound of the communication fault counter \cite{6} to extend the  service time. The obtained delay margin can also be used in designing a delay compensator to improve system performance and stability.
\section{Conclusion}
\label{conclusion}
In this paper, the state space model of a distribution network with PVs considering communication delay has been developed. A delay-dependent time-invariant stability criterion was established by solving the LMI constrains. The delay margins have been obtained for different system parameters. Simulation studies have been carried out on a distribution network with two PVs to verify the calculated delay margins. There are some small differences between the theoretically calculated values and the simulation values due to the errors introduced by the linearization in modeling. The obtained delay margin information is useful in system controller design for a better system performance and a larger stability margin.
\bibliographystyle{IEEEtran}
\bibliography{reference2}

\end{document}